\documentclass[11pt,twoside]{article}
\usepackage{amsmath}
\usepackage{amssymb}
\usepackage{enumerate,epsfig}
\usepackage{epstopdf}
\usepackage{graphicx}
\newtheorem{assumption}{Assumptions}[section]
\newtheorem{theorem}{Theorem}[section]
\newtheorem{lemma}[theorem]{Lemma}

\newtheorem{remark}[theorem]{Remark}
\newtheorem{prop}[theorem]{Proposition}

\newenvironment{proof}{\paragraph{Proof} \phantom{9}}{\hfill$\Box$\bigskip}
\newcommand{\R}{ {\mathbb R} }

\newcommand{\eps}{{\varepsilon}}

\makeatletter
\@addtoreset{equation}{section}
\makeatother

\begin{document}

\title{Selection-Mutation dynamics with spatial dependence}

\author{ Pierre-Emmanuel Jabin,\footnote{CSCAMM and Dept. of Mathematics, University of Maryland,
College Park, MD 20742, USA. P.E. Jabin is partially supported by NSF Grant 1312142 and by NSF Grant RNMS (Ki-Net) 1107444.}
\quad Raymond Strother Schram, \footnote{Dept. of Mathematics, University of Maryland,
College Park, MD 20742, USA.}}

\maketitle

\begin{abstract}
We study the limit of many small mutations for a model of population dynamics. The population is structured by phenological traits and is spatially inhomogeneous. The various sub-populations compete for the same nutrient which diffuses through the spatial environment. 
\end{abstract}
\section{Introduction}
This article solves a conjecture introduced in \cite{MirPer} concerning the limit of many small mutations for some integro-differential model of population dynamics with spatial inhomogeneity.

More precisely, we investigate the behavior as $\eps\rightarrow 0$ of the following model introduced in \cite{LLCEP} and first studied under this scaling in \cite{MirPer}
\begin{align}
&\eps\partial_t n_\eps = (c_\eps r(x) - d(x)\,(1+\varrho_\eps))\,n_\eps +\eps^2\, \partial_{xx} n_\eps(x,y,t),\quad t\geq 0,\ x\in\R,\ y\in\Omega,  \label{eqneps}\\
&\varrho_\eps(y,t)=\int_{\R} n_\eps(x,y,t)\,dx,\label{defrhoeps}\\
&-\Delta_y c_{\eps}(t,y) + [\varrho_{\eps} + \lambda]\,c_{\eps}(t,y) = \lambda\, c_{B},\quad t\geq 0,\ y\in\Omega,\qquad c_\eps|_{\partial\Omega}=0,
\label{eqceps}
\end{align}
where $\Omega$ is a smooth, open and bounded domain of $\R^d$.

In this system $n_\eps(x,y,t)$ represents the density of population with a given trait $x$, at position $y$ in space. The evolution of the population is purely local: Individuals do not move in space but reproduce or die locally. This is modeled in Eq. \eqref{eqneps} by a scaled mutation term $\eps^2\,\partial_{xx} n_\eps$ and a reproduction term $c_\eps r(x) - d(x)\,(1+\varrho_\eps)$. The coefficient $d(x)$ is a death factor, depending on the trait $x$ and combined with the logistic trait $1+\varrho_\eps$ where $\varrho_\eps(y,t)$ is the local total population at the point $y$. This part imposes that the local population cannot grow past a certain limit which is connected to the term $c_\eps\,r(x)$. $r(x)$ is a coefficient giving the reproduction rate in terms of the trait again. $c_\eps(y,t)$ is an external resource (nutrient, oxygen) which the population requires to grow.

The dynamics of each local sub-populations are connected only through this external resource $c_\eps(y,t)$. This nutrient diffuses through space but at a time scales which are assumed to be much faster than the rate of reproduction of the population. Hence the nutrient is assumed to reach its stationary distribution and Eq. \eqref{eqceps} is stationary itself. The nutrient is produced by the environment with a certain flux $\lambda\,c_B$, is used by the population with the absorption term $\varrho_\eps\,c_\eps$ and is also naturally degraded with the additional absorption $\lambda\,c_\eps$.

This model was derived in \cite{LLCEP} to represent the emergence of resistance to drug in cancer therapy, as spatial heterogeneity have been observed and seem important in that context (we refer for example to \cite{AWCQ06}. However the model is also one of the simplest example of coupled dynamics with spatial interaction and selection-mutation. As such it is a good simple framework where to test and illustrate the main properties of such models. 

There are several approaches to study the phenotypical evolution driven by small mutations in replication,  the main objective being to describe the dynamics of the fittest (or dominant) trait
in the population (we refer to \cite{R79} for a general presentation of the biological framework and to \cite{B00} for an early mathematical modeling).

 One of the best known method is probably the so-called adaptive dynamics theory, see \cite{diekmann-04} for an introduction. Adaptive dynamics considers evolution as a succession of invasions by a small mutant population of the population of the dominant trait which is assumed to be at the ecological equilibrium. 

This ansatz can actually be fully justified in some scalings by taking the appropriate limits of stochastic or individually-centered models, see \cite{DL96, CFM06, CFM08, CJM}. Those models are popular in particular since they naturally take fluctuations in the population into account (due to random births, deaths...). Stochastic models can become cumbersome when the total population is too large. 

In that case, it is also possible to derive integro-differential models like the one presented here. This is a slightly different scaling from adaptive dynamics since one does not assume that the population has completely reached the ecological equilibrium in between mutations. Those models naturally take the distribution of traits around the dominant one into account (or can represent multi-modal distributions as in \cite{DBLD07}).
   
The corresponding deterministic method was introduced in \cite{DJMP} and in \cite{DJMR08} without mutations (see also the presentation in \cite{DA14}). This was done in a spatially homogeneous context which already left some open questions. 

In several cases and in particular whenever the competition term is of logistic type as here, then the limiting distribution is concentrated around one trait
\[
n_\eps\longrightarrow \varrho(y,t)\,\delta(x-\bar x(y,t)).
\]
Contrary to the adaptive dynamics scaling, the dominant trait $x(y,t)$ does not solve a closed equation but instead has to calculated from a Hamilton-Jacobi equation with constraint. In the spatially homogeneous case (no spatial dependence on $y$ anywhere), then the rigorous derivation together with the uniqueness at the limit can be found for instance in \cite{BMP09, BarPer, PB07}.

Still in the spatially homogeneous case, when competition is manifested with more complex interactions (see \cite{MPW} for examples of derivation of competitive interactions), then the equilibrium (often called an ESS or evolutionarily stable strategy) may be more complicated. The study of the stability of this ESS can then become an issue of its own. In general the ESS is locally stable, see \cite{R11, R12}, and can be globally stable and unique as in \cite{CJR, JR11} if the competition is strong enough. In this last case, it is possible to still derive the limiting system with the constraint in the Hamilton-Jacobi equation replaced by an implicit formulation, see \cite{ChaJab}. 

Note that the numerical aspects can be interesting and non trivial, see for instance \cite{MPBM11} or \cite{CJL, LCS14}.

As first mentioned above, it has become more and more obvious that spatial heterogeneities are critical for many applications, see \cite{GVA06} for an early example of such of study. And it is natural to try to extend this deterministic theory in those settings, as it is our main goal here.  

The major difficulty in a case such as \eqref{eqneps}-\eqref{eqceps} is that the spatial inhomogeneities could lead to strong oscillations in time. This was the big obstacle in \cite{MirPer}. Of course such oscillations can only develop because the equation on $c_\eps$, \eqref{eqceps}, is stationary in time (elliptic in $y$). If instead it was parabolic in $t$ and $y$ then the difficulty would be mostly avoided thanks to usual regularization in time of parabolic equations (see again \cite{MirPer} or \cite{LMP11}).

Let us now state our main result. Following \cite{DJMP}, we introduce
\begin{equation}
u_\eps=\eps\,\log n_\eps,\quad\mbox{or}\ n_\eps=e^{u_\eps/\eps}.\label{defueps}
\end{equation}
From Eq. \eqref{eqneps}, we see that $u_\eps$ satisfies
\begin{equation}
\partial_t u_\eps=c_\eps r(x) - d(x)\,(1+\varrho_\eps)+\eps\,\partial_{xx} u_\eps+|\partial_x u_\eps|^2.\label{equeps}
\end{equation} 
\begin{assumption}
Assume that $r,\;d\in W^{2,\infty}(\R)$ and are bounded from below and from above. Define
\[
\overline \varrho=c_B\,\frac{\sup r}{\inf d}-1,\quad \underline \varrho=\frac{\inf r}{\sup d}\,\frac{c_B\,\lambda}{\lambda + \overline\varrho}-1,
\] 
and assume that $c_B>0$ is large enough so that $0<\underline\varrho<\overline\varrho$. Assume moreover that $r$ is \emph{concave} and that $d$ is \emph{convex}.
\label{assumpt1}
\end{assumption}
We can prove
\begin{theorem} Assume \eqref{assumpt1}, that $u_\eps^0$ is uniformly Lipschitz in $x$ and is uniformly strictly concave in $x$ namely
\[
\sup_\eps\sup_{y,x} |\partial_x u_\eps^0(x,y)|<\infty,\quad -\frac{1}{K^0}\leq \partial_{xx} u^0_\eps(x,y)\leq -K^0<0,
\]
for some constant $K>0$ independent of $\eps$. Assume that $n^0_\eps\rightarrow \varrho^0(y)\,\delta(x-x^0(y))$ for some $x^0(y)$. Define $\bar\varrho^0,\;\bar c^0$ as the unique solution (see section \ref{secmetastable}) to  
\[
\begin{split}
&-\Delta_y\bar c^0 + (\lambda +\bar{\varrho}^0)\,\bar c^0 = \lambda\, c_B,\\
&\bar c^0\, r(x^0(y)) - (\bar \varrho^0 + 1)\,d(x^0(y))=0.
\end{split}
\]
There exists $\Lambda$ s.t. if 
\[
\int_\Omega |\varrho^0-\bar\varrho^0|^2\,dy<\Lambda,
\]
then up to an extracted subsequence, for any $T>0$, $c_\eps$ converges to $c$ in $L^2([0,\ T]\times\Omega)$, $\int_\Omega \varrho_\eps\,\varphi\,dy$ converges to $\int_\Omega\varrho\,\varphi\,dy$ in $C([0,\ T])$; $u_\eps$ converges weak-* in $L^\infty([0,\ T]\times\Omega\times\R)$ to $u$ which is a viscosity solution to
\[
\partial_t u=c(y,t)\,r(x)-d(x)\,(1+\varrho(y,t))+|\partial_x u|^2,\quad t\geq0,\ x\in\R,\ y\in\Omega, 
\] 
with the constraint $\max_x u(.,y,t)=0$ and coupled to 
\[
\Delta_y c(y,t)+[\varrho(y,t)+\lambda]\,c(y,t)=\lambda\,c_B,\quad t\geq 0,\ y\in\Omega,\quad c|_{\partial\Omega}=0.
\]
\label{mainresult}
\end{theorem}
\begin{remark}
We stated the limiting equation as a Hamilton-Jacobi equation with the constraint $\max u=0$. However we could also have formulated in terms of the meta-stable states defined in section \ref{secmetastable}. In particular the system \ref{metastable} would let us calculate directly $\varrho$ and $c$ in terms of the maximum point of $u$. Those meta-stable states and the corresponding two-time scales analysis are at the heart of the proof and were inspired by the approach in \cite{ChaJab, Jab}.
\end{remark}
\begin{remark}
The reason for the assumption $\int_\Omega|\varrho^0-\bar\varrho|^2\,dy<\Lambda$ is the analysis of the stability of the meta-stable states. Unfortunately, we are only able to prove their local stability and this is what imposes this smallness assumption.
\end{remark}
\begin{remark}
The same type of ideas as presented here can be applied to several variants of the system \eqref{eqneps}-\eqref{eqceps}. The mutations for instance can be replaced by an integral kernel without substantial modifications in the proof.
\end{remark}
\section{A Priori Estimates}
The first step is to find bounds on $c_\eps$ and $\varrho_\eps$. Those mostly follow the a priori estimates developed in \cite{MirPer}.  We first draw some basic estimates that show these are non-negative and bounded, which we then use to form more precise estimates.

We start with some straightforward time uniform bounds
\begin{prop}
Under the assumptions \ref{assumpt1}, one has the following bounds on the solution to System \eqref{eqneps}-\eqref{eqceps}, 
\[
\underline\varrho\leq \varrho_\eps\leq \overline\varrho ,\qquad \frac{c_B\,\lambda}{\lambda + \overline\varrho} \leq c_\eps\leq c_B.
\]
Moreover $c_\eps$ is uniformly Lipschitz and in $W^{2,p}_{loc}$ for any $p$, {\em i.e.} for any $p<\infty$ and any compact set $\Omega$
\[
\sup_{\eps,t} \left(\|\nabla_y c_\eps\|_{L^\infty_y}+\|\nabla_y c_\eps\|_{L^p_y(\Omega)}\right)<\infty.  
\]
\label{aprioritimeunif}
\end{prop}
\begin{proof}
Since $c_B>0$, by the maximum principle on Eq. \eqref{eqceps}, one has immediately that $c_\eps\geq 0$.
 
Now integrate over $x$ (the space of traits) Eq. \eqref{eqneps}, using the bounds on $r$ and $d$
\begin{equation}
\eps\partial_t \varrho_\eps \geq \inf r\, c_\eps\,\varrho_\eps - \sup d\,\varrho_\eps(\varrho_\eps +1). \label{intermediary1}
\end{equation}
This implies that $\varrho_\eps\geq 0$.

Coming back to Eq. \eqref{eqceps}, one may use the maximum principle again  to find that
\[
\sup [\varrho_\eps + \lambda]\,c_\eps \leq \lambda\, c_B,
\]
and together with the non-negativity of $\varrho_\eps$, $c_\eps \leq c_B$. 

Integrating again Eq. \eqref{eqneps} in $x$ and using now this upper bound on $c_\eps$, 
\[
\eps\partial_t \varrho_\eps \leq \sup r\, c_B\,\varrho_\eps - \inf d_\eps\, (1+\varrho_\eps)\,\varrho_\eps. 
\]
Observe that the right-hand side becomes negative if $\varrho_\eps>\overline\varrho$. Thus $\varrho_\eps\leq \overline\varrho$. 

This in turn enables to obtain a more precise lower bound on $c_\eps$. Using the maximum principle again on \eqref{eqceps}
\[
 \lambda c_B\leq  [\overline\varrho + \lambda]\,\inf c_\eps,
\]
leading to $c_\eps\geq \lambda\,c_B/(\lambda+\overline\varrho)$. Introducing this bound in Eq. \eqref{intermediary1},
\[
\eps\partial_t \varrho_\eps \geq \inf r\, \frac{\lambda\,c_B}{\lambda+\overline\varrho}\,\varrho_\eps - \sup d\,\varrho_\eps(\varrho_\eps +1),
\]
and the right-hand side is positive if $\varrho_\eps>\underline\varrho$ which finishes the proof. Note that it would be possible to improve the upper bound on $c_\eps$ by using this improved lower bound on $\varrho_\eps$ giving 
\[
 c_\eps \leq \frac{c_B\lambda}{\lambda + \underline\varrho}.
\]
Finally introducing those bounds in Eq. \eqref{eqceps} shows that $\Delta_y c_\eps$ is bounded which concludes by standard elliptic estimates.
\end{proof}

We now turn to a second set of {\em a priori} estimates on $u_\eps$ defined through and which are non uniform in time.
\begin{prop}
Assume \eqref{assumpt1}, that $u_\eps^0$ is uniformly Lipschitz in $x$ and is uniformly strictly concave in $x$ namely
\[
\sup_\eps\sup_{y,x} |\partial_x u_\eps^0(x,y)|<\infty,\quad -\frac{1}{K^0}\leq \partial_{xx} u^0_\eps(x,y)\leq -K^0<0,
\]
for some constant $K>0$ independent of $\eps$. 
Then this is preserved for any finite time and $u_\eps$ is Lipschitz in time:
\[
\sup_{t\in[0,\ T]}\sup_{\eps,y,x} (|\partial_x u_\eps(x,y,t)|+|\partial_t u_\eps(x,y,t)|)<\infty,\quad -\frac{1}{K(t)}\leq \partial_{xx} u_\eps(x,y)\leq -K(t)<0,
\]
for any $T>0$ and for some  $K(t)<0$ for any $t$ and independent of $\eps$. Moreover one has the following estimate on the maximum in $x$ of $u_\eps$
\[
-\frac{\eps}{2}\, \log\eps+\eps\log \frac{\underline\varrho\,\sqrt{K}}{\sqrt{2\,\pi}}\leq \max_x u_\eps(.,y,t)\leq -\frac{\eps}{2}\, \log\eps+\eps\log \frac{\overline\varrho}{\sqrt{2\,K\,\pi}}.
\]\label{aprioriueps}
\end{prop}
\begin{proof}
First, let us differentiate once in $x$ Eq. \eqref{equeps}
\[
\partial_t \partial_{x} u_\eps = c_\eps\, r^{\prime} - d^{\prime}\,(\varrho_\eps +1)
+ 2\partial_x u_\eps\,\partial_{xx}u_\eps + 
\eps \partial_{xx} \partial_{x}u_\eps.
\]
By the maximum principle, one has that
\[
\partial_t \sup_{x,y} |\partial_x u_\eps|\leq \sup_{x,y} (c_\eps\,|r'|+|d'|\,(\varrho_\eps+1).
\]
Using now the bounds from Prop. \ref{aprioritimeunif} and the Lipschitz bounds on $r$ and $d$
\[
\partial_t \sup_{x,y} |\partial_x u_\eps|\leq C\,(c_B+\overline\varrho+1), 
\] 
for some $C>0$ which gives the uniform Lipschitz bound in $x$.

Differentiate now twice Eq. \eqref{equeps} to find
\[
\partial_t \partial_{xx} u_\eps = c_\eps r^{\prime\prime} - d^{\prime\prime}(\varrho_\eps +1)
+ 2|\partial_{xx}u_\eps|^2 + 2\partial_x u_\eps \partial_x \partial_{xx} u_\eps + 
\eps \partial_{xx} \partial_{xx}u_\eps.
\]
By assumption \eqref{assumpt1} (more precisely the concavity of $r$ and the convexity of $d$) and since $c_\eps>0$ and $\varrho_\eps>0$, this leads to
\[
\partial_t \partial_{xx} u_\eps \leq 2|\partial_{xx}u_\eps|^2 + 2\partial_x u_\eps \partial_x \partial_{xx} u_\eps + 
\eps \partial_{xx} \partial_{xx}u_\eps.
\]
The inequality contains a Riccati-like term, $|\partial_{xx}u_\eps|^2$, and may in general lead to blow-up. However we have assumed that $\partial_{xx} u_\eps^0<0$ which precisely guarantees that blow-up do not occur. In fact one immediately has that $\partial_{xx} u_\eps(x,y,t)\leq 0$ for every $t>0$, and $x,\;y$. But by the maximum principle, one has the more precise
\[
\partial_t \sup_{x,y} \partial_{xx} u_\eps\leq 2|\sup_{x,y} \partial_{xx} u_\eps|^2,
\]
and thus $\sup_\eps\sup_{x,y}\partial_{xx} u_\eps\leq -K(t)$ for some function $K(t)$ which may converge to $0$ as $t\rightarrow \infty$ but remains strictly positive in the meantime.

The lower bound on $\partial_{xx} u_\eps$ is even easier to obtain, just by using that $r''$ and $d''$ are bounded, again from assumption \eqref{assumpt1}.

Then from the bounds on $\partial_x u_\eps$, $\partial_{xx} u_\eps$, Eq. \eqref{equeps} directly provides a uniform bound on $\partial_t u_\eps$.

It only remains to obtain the bounds on the maximum of $u_\eps$.  Observe that from the uniform bound on $\partial_{xx} u_\eps$, for any $t$ and $y$, $u_\eps$ necessarily has a unique maximum in $x$. 
Denote
\[
x_\eps(y,t)=\mbox{argmax}\, u_\eps(.,y,t).
\]
One has that
\[
\varrho_\eps = \int_\R n_\eps = \int_\R e^{\frac{u_\eps(x,y,t)}{\eps}}\,dx. 
\]
But by Taylor-Lagrange, $u_\eps(x,y,t)=u_\eps(x_\eps(y,t),y,t)-\frac{1}{2}\partial_{xx} u_\eps(z,y,t)\,(x-x_\eps)^2$. Thus
\[
u_\eps(x_\eps(y,t),y,t)-\frac{1}{2\,K}\,(x-x_\eps)^2\leq u_\eps(x,y,t)\leq u_\eps(x_\eps(y,t),y,t)-\frac{K}{2}\,(x-x_\eps)^2,
\]
and
\[
\int_\R e^{\frac{u_\eps( x_\eps,y,t) - \frac{1}{2\,K}(x-x_\eps)^2}{\eps}}\,dx\leq \varrho_\eps\leq \int_\R e^{\frac{u_\eps(x_\eps,y,t) - \frac{K}{2}(x-x_\eps)^2}{\eps}}\,dx. 
\]
Therefore
\[
e^{\max_x u_\eps(.,y,t)/\eps}\,\sqrt{2\,K\,\eps\,\pi}\leq \varrho_\eps\leq e^{\max_x u_\eps(.,y,t)/\eps}\,\frac{\sqrt{2\,\,\eps\,\pi}}{K},
\]
which concludes the proof by using $\underline\varrho\leq\varrho_\eps\leq \overline\varrho$ from Prop. \ref{aprioritimeunif}.
\end{proof}

Prop. \ref{aprioritimeunif} and Prop. \ref{aprioriueps} provide strong bounds but no compactness in time on $\varrho_\eps$ or $c_\eps$. As explained in the introduction, this is the key problem as it prevents us from passing to the limit in the product $c_\eps\,\varrho_\eps$ in Eq. \eqref{eqceps} for instance.
\section{The meta-stable state\label{secmetastable}}
One solution to the question of time oscillations is to introduce the right meta-stable states to the system and compare them to the solution. This is what we are doing in this section.
\subsection{Time regularity of the maximum point of $u_\eps$}
We recall that $x_\eps(y,t)$ is defined as the unique point where $u_\eps(.,y,t)$ attains its maximum.

From the regularity of $u_\eps$ provided by Prop. \ref{aprioriueps}, it is possible to deduce the following time regularity of $x_\eps$
\begin{lemma}
Assume \eqref{assumpt1} and all the assumptions of Prop. \ref{aprioriueps}, then for any $T>0$, there exists a constant $C_T>0$ uniform in $\eps$ s.t.
\[
|x_\eps(t,y)-x_\eps(s,y)|\leq C_T\,(\sqrt{\eps}+\sqrt{|t-s|}),\quad\forall t,\;s\in[0,\ T].
\]\label{regxeps}
\end{lemma}
\begin{proof}
Observe that by Prop. \ref{aprioriueps}, for some constant $C_T^1$, depending only on $K(t)$ on $[0,\ T]$ and hence independent of $\eps$
\[
|u_\eps(x_\eps(y,t),y,t)-u_\eps(x_\eps(y,s),y,s)|\leq C_T^1\,\eps,
\]
for any $t,\;s$ in $[0,\ T]$. On the other hand by the Lipschitz in $t$ bound
\[
|u_\eps(x_\eps(y,t),y,t)-u_\eps(x_\eps(y,t),y,s)|\leq C_T^2\,|t-s|,
\]
and from the upper bound on $\partial_{xx} u_\eps$
\[
u_\eps(x_\eps(y,s),y,s)-u_\eps(x_\eps(y,t),y,s)\geq \frac{K(T)}{2}\,|x_\eps(y,s)-x_\eps(y,t)|^2.
\]
Therefore
\[
\frac{K(T)}{2}\,|x_\eps(y,s)-x_\eps(y,t)|^2\leq C_T^1\,\eps+ C_T^2\,|t-s|,
\]
finishing the proof.
\end{proof}
\subsection{Definition and basic properties of the meta-stable states}
We define the meta-stable states $(\bar\varrho_\eps(y,t),\;\bar c_\eps(y,t))$ through the following system
\begin{equation}\begin{split}
&-\Delta_y\bar c_\eps + (\lambda +\bar{\varrho}_\eps)\,\bar c_\eps = \lambda\, c_B,\\
&\bar c_\eps r(x_\eps(y,t)) - (\bar \varrho_\eps + 1)\,d(x_\eps(y,t))=0.
\end{split}\label{metastable}
\end{equation} 
Note that the two equations are in fact stationary in time. The time dependence in $\bar\varrho_\eps$ and $c_\eps$ is only due to the fact that the functions $r$ and $d$ in \eqref{metastable} is taken at the point $x_\eps(y,t)$. 

From the second equation in \eqref{metastable}, it is possible to reduce the system to a non-linear elliptic equation. Observe indeed that defining
\[
f_\eps(y,t)=\frac{r(x_\eps(y,t))}{d(x_\eps(y,t))},
\]
one has obviously
\begin{equation}
\bar\varrho_\eps=\bar c_\eps\,f_\eps-1>0,\label{algrhoc}
\end{equation}
from assumption \eqref{assumpt1}.

Plugging this into the equation on $\bar c_\eps$ yields
\begin{equation}
-\Delta_y \bar c_\eps + (\lambda + \bar c_\eps\, f_\eps - 1)\,\bar c_\eps = \lambda c_B.\label{ellipticceps}
\end{equation}
It is straightforward to obtain
\begin{lemma}
There exists a unique solution to the System \eqref{metastable} with $c_\eps\in L^\infty(\R_+,\ L^\infty(\Omega)\cap H^1(\Omega))$ and $\bar\varrho_\eps\in L^\infty(\R_+\times \Omega)$. In addition those solutions satisfy
\[
\underline \varrho\leq \bar \varrho_\eps\leq \overline\varrho,\quad \frac{C_B\,\lambda}{\lambda+\overline\varrho}\leq \bar c_\eps\leq c_B.
\]\label{existmetastable}
\end{lemma} 
Note that time is only a parameter in Eq. \eqref{ellipticceps} and thus the proof of Lemma \ref{existmetastable} is a straightforward application of standard non-linear elliptic techniques. The specific bounds on $\bar c_\eps$ and $\bar\varrho_\eps$ follow exactly the lines of the proof of Prop. \ref{aprioritimeunif}. For this reason we skip the proof and instead present a variant of the uniqueness argument which provides regularity in time.
\begin{prop}
Assume \eqref{assumpt1} and all the assumptions of Prop. \ref{aprioriueps}, then for any $T>0$, there exists a constant $\bar C_T>0$ uniform in $\eps$ s.t. for any $t,\;s\in[0,\ T]$
\[
\begin{split}
&\|\bar c_\eps(.,t)-\bar c_\eps(.,s)\|_{L^\infty_y}\leq \bar C_T\,(\sqrt{\eps}+\sqrt{|t-s|}),\\
&\|\bar \varrho_\eps(.,t)-\bar \varrho_\eps(.,s)\|_{L^\infty_y}\leq \bar C_T\,(\sqrt{\eps}+\sqrt{|t-s|}).
\end{split}
\]\label{regrhoeps}
\end{prop}
\begin{proof}
As suggested above the proof revolves around uniqueness estimates for Eq. \eqref{ellipticceps}. Denote
\[
\delta \bar c_\eps=\bar c_\eps(y,t)-\bar c_\eps(y,s).
\]
Then from Eq. \eqref{ellipticceps}
\[
-\Delta_y \delta\bar c_\eps + (\lambda + (\bar c_\eps(y,t)+\bar c_\eps(y,s))\, f_\eps(y,t) - 1)\,\delta \bar c_\eps =\bar c_\eps^2(y,s)\,(f_\eps(y,s)-f_\eps(y,t)).
\]
On the other hand since $r$ and $d$ are Lipschitz, one has by Lemma \ref{regxeps}
\begin{equation}
|f_\eps(y,s)-f_\eps(y,t)|\leq C_T\,(\sqrt{\eps}+\sqrt{|t-s|}),\label{fepsdiff}
\end{equation}
and hence
\[
\left|-\Delta_y \delta\bar c_\eps + (\lambda + (\bar c_\eps(y,t)+\bar c_\eps(y,s))\, f_\eps(y,t) - 1)\,\delta \bar c_\eps\right|\leq  c_B^2\,C_T\,(\sqrt{\eps}+\sqrt{|t-s|}).
\]
The coefficient $\lambda + (\bar c_\eps(y,t)+\bar c_\eps(y,s))\, f_\eps(y,t) - 1)$ is bounded from below thanks to the a priori estimates in Lemma \ref{existmetastable}
\[
\lambda + (\bar c_\eps(y,t)+\bar c_\eps(y,s))\, f_\eps(y,t)-1\geq \lambda+1+2\,\underline\varrho.
\]
Standard elliptic estimates thus imply that
\[
\|\delta\bar c_\eps\|_{W^{2,p}_{loc}}\leq C_{T,p} \,(\sqrt{\eps}+\sqrt{|t-s|}),
\]
for any $p<\infty$ and for some constant $C_{T,p}>0$. This concludes the bound on $\delta\bar c_\eps$. The bound on $\bar\varrho_\eps(y,t)-\bar\varrho_\eps(y,s)$ follows directly from formula \eqref{algrhoc} and the regularity property \eqref{fepsdiff} for $f_\eps$.
\end{proof}
%
\section{Comparison with the meta-stable states}
\subsection{An approximate equation on $\varrho_\eps$}
Now we will compare $\bar\varrho_\eps$ and $\bar c_\eps$ with $\varrho_\eps$ and $c_\eps$. Ideally one would have a closed system of equations on $\varrho_\eps$ and $c_\eps$ but this is of course not possible in general. Nevertheless it is possible to use the concentration in $n_\eps$ to derive an approximate closed system.

Let us begin with this concentration property
\begin{lemma}
Assume \eqref{assumpt1} and the assumptions of Prop. \ref{aprioriueps}. Then for any $T$, there exists $C_T>0$ s.t. for any $\varphi\in W^{1,\infty}(\R)$ and any $t\in[0,\ T]$, $y\in\Omega$
\[
\left|\int_\R \varphi(x)\,n_\eps(x,y,t)\,dx-\varphi(x_\eps(y,t))\,\varrho_\eps(y,t)\right| \leq 
C_T\,\sqrt{\eps}\, \|\varphi\|_{W^{1,\infty}}\,.
\]\label{concentration}
\end{lemma}
\begin{proof}
Of course
\[
|\varphi(x)-\varphi(x_\eps(y,t))|\leq \|\varphi\|_{W^{1,\infty}}\,|x-x_\eps(y,t)|,
\]
and
\[\begin{split}
&\int_\R |x-x_\eps|\,n_\eps(x,y,t)\,dx=\int_\R |x-x_\eps|\,e^{u_\eps(x,y,t)/\eps}\,dx\\
&\leq e^{u_\eps(x_\eps,y,t)/\eps}\,\int_\R |x-x_\eps|\,e^{-K(t)\,|x-x_\eps|^2/2\eps}\,dx,
\end{split}
\]
by the concavity of $u_\eps$ from Prop. \ref{aprioriueps}. Therefore using the estimate on the maximum of $u_\eps$ also from Prop. \ref{aprioriueps}
\[
\int_\R |x-x_\eps|\,n_\eps(x,y,t)\,dx\leq \frac{\overline\varrho}{\sqrt{2K}}\,\sqrt\eps\,\int_\R |z|\,e^{-|z|^2}\,dz, 
\]
which proves the lemma.
\end{proof}
 
It is now straightforward to derive the equation on $\varrho_\eps$ 
\begin{prop}
Assume \eqref{assumpt1} and the assumptions of Prop. \ref{aprioriueps}. Then there exists $R_\eps(y,t,s)$ and $C_T>0$ independent of $\eps$ with 
\[
|R_\eps(y,t,s)|\leq C_T\,(\sqrt{\eps}+\sqrt{|t-s|}),\quad\forall y\in\Omega,\quad \forall t,\;s\in [0,\ T],
\] 
and such that
\begin{equation}
\eps\,\partial_t\varrho_\eps(t,y) = (r(x_\eps(y,s))\,c_\eps(t,y) - d(x_\eps(y,s))\,(1+\varrho_\eps))\,\varrho_\eps+R_\eps(y,t,s).\label{eqrhoeps}
\end{equation}\label{propeqrhoeps}
\end{prop}
\begin{proof}
Start by integrating Eq. \eqref{eqneps} in $x$
\[
\eps\,\partial_t\varrho_\eps(t,y) = \int[r(x)\,c_\eps(t,y) - d(x)\,(1+\varrho_\eps)]\,n_\eps\,dx.
\]
Simply by using Lemma \ref{concentration} and the Lipschitz bound on $r$ and $d$ we have that
\[
\eps\,\partial_t\varrho_\eps(t,y) = (r(x_\eps(y,t))\,c_\eps(t,y) - d(x_\eps(y,t))\,(1+\varrho_\eps))\,\varrho_\eps+R_{\eps,1}(y,t),
\]
with 
\[
|R_{\eps,1}(y,t)|\leq C_T\,\sqrt{\eps}.
\]
It now only remains to change $r(x_\eps(y,t))$ in $r(x_\eps(y,s))$ and similarly for $d$ to derive Eq. \eqref{eqrhoeps}. This is again due to the Lipschitz bound on $r$ and $d$ and the time regularity on $x_\eps$ provided by Lemma \ref{regxeps}.
\end{proof}

Observe that in Eq. \eqref{eqrhoeps} we have delocalized in time the coefficients in $r$ and $d$. This is actually a critical step as it will allow to use a simple stability estimate for the corresponding dynamical system.
\subsection{Stability of Eqs. \eqref{eqrhoeps}-\eqref{eqceps}}
The above derivation naturally leads to the study in large time of the given system
\begin{equation}
\begin{split}
&\partial_t\varrho(t,y) = (r(x_\eps(y,s))\,c(t,y) - d(x_\eps(y,s))\,(1+\varrho))\,\varrho,\\
&-\Delta_y c(t,y) + [\varrho + \lambda]\,c(t,y) = \lambda\, c_{B},
\end{split}\label{odesystem}
\end{equation}
for a fixed $s$. As one can see from their definition \eqref{metastable}, the meta-stable states $\bar\varrho_\eps,\;\bar c_\eps$ are fixed point of this system. The key in our analysis is that they are in fact stable fixed points, at least linearly ({\em i.e.} close to equilibrium). This can be seen very simply through a Lyapunov functional
\begin{lemma}
Assume \eqref{assumpt1} and the assumptions of Prop. \ref{aprioriueps}. There exists a constant $C>0$ s.t any solution $\varrho,\;c$ to System \eqref{odesystem} with $\underline\varrho\leq\varrho\leq\overline\varrho$ satisfies
\[\begin{split}
&\frac{d}{dt}\frac{1}{2}\int|\varrho(y,t)-\bar{\varrho_\eps}(y,s)|^2\, \frac{\bar c_\eps(y,s)}{\bar{\varrho_\eps}(y,s)\,r(x_\eps(y,s))}\,dy\leq-C\,\int |\varrho(y,t)-\bar{\varrho_\eps}(y,s)|^2\,dy\\
&\qquad-\int|\nabla(c(y,t)-\bar{c_\eps}(y,s))|^2\,dy-C\,\int (c(y,t)-\bar{c_\eps}(y,s))^2\,dy\\
&\qquad
+\int (\varrho(y,t)- \bar{\varrho_\eps}(y,s))^2\, (c(y,t)-\bar{c_\eps}(y,s))\,\frac{\bar{c_\eps}(y,s)}{\bar{\varrho_\eps}(y,s)}\,dy.
\end{split}
\]\label{lyapunov}
\end{lemma} 
Remark that this is a Lyapunov functional only when $\varrho,\; c$ is close to $\bar\varrho_\eps,\;\bar c_\eps$ as then the only term which can be positive in the right-hand side, $\int (\varrho- \bar{\varrho_\eps})^2\, (c-\bar{c_\eps})\,\frac{\bar{c_\eps}}{\bar{\varrho_\eps}}\,dy$ is negligible in front of the others. We will of course heavily use this property later on.

\begin{proof}
It is a straightforward calculation from th equations \eqref{odesystem}.
In the following to avoid writing all the variables at each step, $x_\eps$, $\bar\varrho_\eps$ and $\bar c_\eps$ are taken at the point $(y,s)$ while $\varrho$ and $c$ are always taken at the point $(y,t)$. Using the first equation
\[
\frac{d}{dt}\frac{1}{2}\int|\varrho-\bar{\varrho_\eps}|^2 \frac{\bar c_\eps}{\bar{\varrho_\eps}\,r(x_\eps)}\,dy=\int[\varrho_\eps - \bar{\varrho_\eps}]\,[r(x_\eps)\,c -d(x_\eps)(1+\varrho)]\, \frac{\varrho\,\bar c_\eps} {\bar{\varrho_\eps}\,r(x_\eps)}\,dy.
\]
Recall that
\[
r(x_\eps)\,\bar c_\eps -d(x_\eps)(1+\bar\varrho_\eps)=0,
\]
so that
\[\begin{split}
\frac{d}{dt}\frac{1}{2}\int|\varrho-\bar{\varrho_\eps}|^2 \frac{\bar c_\eps}{\bar{\varrho_\eps}\,r(x_\eps)}\,dy&=\int[\varrho_\eps - \bar{\varrho_\eps}]\,[r(x_\eps)\,(c-\bar c_\eps) -d(x_\eps)(\varrho-\varrho_\eps)]\, \frac{\varrho\,\bar c_\eps} {\bar{\varrho_\eps}\,r(x_\eps)}\,dy\\
&=-\int |\varrho-\bar{\varrho_\eps}|^2\frac{\bar{c_\eps}\,d(x_\eps)\, \varrho} {\bar{\varrho_\eps}\,r(x_\eps)} + \int(\varrho- \bar{\varrho_\eps})\, (c-\bar{c_\eps})\,\frac{\bar{c_\eps}\,\varrho}{\bar{\varrho_\eps}}.
\end{split}
\]
Write simply 
\[
 \int(\varrho- \bar{\varrho_\eps})\, (c-\bar{c_\eps})\,\frac{\bar{c_\eps}\,\varrho}{\bar{\varrho_\eps}}\,dy=\int(\varrho- \bar{\varrho_\eps})\,(c-\bar{c_\eps})\,\bar c_\eps\, dy+\int (\varrho- \bar{\varrho_\eps})^2\, (c-\bar{c_\eps})\,\frac{\bar{c_\eps}}{\bar{\varrho_\eps}}\,dy.
\]
Use the second equations of \eqref{odesystem} and \eqref{metastable} to find
\[
-\Delta(c - \bar{c_\eps}) + \lambda\,(c - \bar{c_\eps}) + \varrho c - \bar{\varrho_\eps}\bar{c_\eps} = 0
\]
Derive the energy estimate by multiplying by $c-\bar c_\eps$ and integrating in $y$
\[
\int|\nabla(c-\bar{c_\eps})|^2\,dy + \lambda\,\int(c-\bar{c_\eps})^2\,dy + \int(c-\bar{c_\eps})(\varrho\, c - \bar{\varrho_\eps}\bar{c_\eps}) = 0.
\]
Of course
\[
\int(c-\bar{c_\eps})(\varrho\, c - \bar{\varrho_\eps}\bar{c_\eps})\,dy=\int\bar{c_\eps}\,(c-\bar{c_\eps})\,(\varrho -\bar{\varrho_\eps})\,dy+\int\varrho\,(c_\eps -\bar{c_\eps})^2\,dy.
\]
Therefore
\[
\int(\varrho- \bar{\varrho_\eps})\,(c-\bar{c_\eps})\,\bar c_\eps\, dy=\int|\nabla(c-\bar{c_\eps})|^2\,dy + \int(\lambda+\varrho)\,(c-\bar{c_\eps})^2\,dy.
\]
Introducing in the first estimate yields
\[\begin{split}
\frac{d}{dt}\frac{1}{2}\int|\varrho-\bar{\varrho_\eps}|^2 \frac{\bar c_\eps}{\bar{\varrho_\eps}\,r(x_\eps)}\,dy&=-\int |\varrho-\bar{\varrho_\eps}|^2\frac{\bar{c_\eps}\,d(x_\eps)\, \varrho} {\bar{\varrho_\eps}\,r(x_\eps)}-\int|\nabla(c-\bar{c_\eps})|^2\,dy\\
&-\int(\lambda+\varrho)\,(c-\bar{c_\eps})^2\,dy
+\int (\varrho- \bar{\varrho_\eps})^2\, (c-\bar{c_\eps})\,\frac{\bar{c_\eps}}{\bar{\varrho_\eps}}\,dy.
\end{split}
\]
To conclude the proof, it is enough to use the upper and lower bounds on $\varrho,\;\bar\varrho_\eps,\;\bar c_\eps$.
\end{proof}
\subsection{The final estimate between $(\varrho_\eps,\;c_\eps)$ and $(\bar \varrho_\eps,\;\bar c_\eps)$}
We are now ready to estimate the difference between $(\varrho_\eps,\;c_\eps)$ and $(\bar \varrho_\eps,\;\bar c_\eps)$ as per
\begin{prop}
Assume \eqref{assumpt1} and the assumptions of Prop. \ref{aprioriueps}. There exists $\Lambda$ s.t. if at $t=0$
\[
\int_{\Omega} |\varrho_\eps(y,t=0)-\bar\varrho_\eps(y,t=0)|^2\,dy<\Lambda,
\]
then for any $T>0$, there exists $\eps_T$ s.t. for any $\eps<\eps_T$ the same is true for all times $t\in[0,\ T]$ and moreover for some $C_T>0$ independent of $\eps$ and for any $t\in [\sqrt{\eps},\ T]$
\[
\int_{\Omega} |\varrho_\eps(y,t)-\bar\varrho_\eps(y,t)|^2\,dy+\int_0^T \|c_\eps(y,s)-\bar c_\eps(y,s)\|_{H^1(\Omega)}^2\,ds\leq C_T\,\sqrt\eps.
\] \label{propcompare}
\end{prop}
Remark that $\eps^{1/2}$ is chosen just for convenience, any $\eps^{k}$ with $k<1$ would work. Remark as well that the last estimates holds only for $t\geq \eps^{1/2}$; it is indeed due to dissipation properties of the system which need a minimal time to be effective.
\begin{proof}
We define
\[
M_\eps=\frac{\bar c_\eps(y,s)}{\bar{\varrho_\eps}(y,s)\,r(x_\eps(y,s))}.
\]
We start by using Lemma \ref{lyapunov} but of course on $\varrho_\eps$. Given that $\varrho_\eps$ solves Eq. \eqref{eqrhoeps} instead of the first equation of \eqref{odesystem}, one has that
\[\begin{split}
&\eps\,\frac{d}{dt}\frac{1}{2}\int|\varrho_\eps(y,t)-\bar{\varrho_\eps}(y,s)|^2\, M_\eps(s,y)\,dy\leq-C\,\int |\varrho_\eps(y,t)-\bar{\varrho_\eps}(y,s)|^2\,dy\\
&\qquad-\int|\nabla(c_\eps(y,t)-\bar{c_\eps}(y,s))|^2\,dy-C\,\int (c_\eps(y,t)-\bar{c_\eps}(y,s))^2\,dy\\
&\qquad
+\int (\varrho_\eps(y,t)- \bar{\varrho_\eps}(y,s))^2\, (c_\eps(y,t)-\bar{c_\eps}(y,s))\,\frac{\bar{c_\eps}(y,s)}{\bar{\varrho_\eps}(y,s)}\,dy+\bar R_\eps(t,s),
\end{split}
\]
with
\[
\bar R_\eps(t,s)=\int_{\Omega} (\varrho_\eps(y,t)-\bar{\varrho_\eps}(y,s))\,R_\eps(y,t,s)\,dy.
\]
By Lemma \ref{propeqrhoeps},
\[
\bar R_\eps\leq C_T\,(\sqrt{\eps}+\sqrt{|t-s|})\,\|\varrho_\eps-\bar\varrho_\eps\|_{L^2_y}\leq \tilde C_T\,(\eps+|t-s|)+\frac{C}{2}\int |\varrho_\eps(y,t)-\bar\varrho_\eps(y,s)|^2\,dy.
\]
Let us now handle the term of order $3$. By H\"older estimate for any $p>\geq1$
\[\begin{split}
\int(\varrho_\eps - \bar{\varrho_\eps})^2\, (c_\eps-\bar{c_\eps})\,\frac{\bar c_\eps}{\bar\varrho_\eps}\,dy &\leq \frac{c_B}{\underline\varrho}\,\left(\int (\varrho_\eps - \bar{\varrho_\eps})^2\,|c_\eps-\bar{c_\eps}|^p\,dy\right)^{1/p}\,\left( \int|\varrho_\eps - \bar{\varrho_\eps}|^2\,dy\right)^{1-1/p}\\
&\leq \frac{c_B}{\underline\varrho}\,(2\,\overline\varrho)^{1/p} \|c_\eps-\bar{c_\eps}\|_{L^p_y}\,\left( \int|\varrho_\eps - \bar{\varrho_\eps}|^2\,dy\right)^{1-1/p}.
\end{split}
\]
Choose now $p>2$ s.t. by Sobolev embedding
\[
\|c_\eps-\bar{c_\eps}\|_{L^p_y}\leq C_d\,\|c_\eps-\bar{c_\eps}\|_{H^1_y}.
\]
Consequently for some constant $\bar C$ depending only on the dimension, $c_B$, $\overline\varrho$ and $\underline\varrho$
\[\begin{split}
\int(\varrho_\eps - \bar{\varrho_\eps})^2\, (c_\eps-\bar{c_\eps})\,\frac{\bar c_\eps}{\bar\varrho_\eps}\,dy &\leq \frac{C}{2}\,\|c_\eps-\bar{c_\eps}\|_{H^1_y}+\bar C\,\left( \int|\varrho_\eps - \bar{\varrho_\eps}|^2\,dy\right)^{2-2/p}.
\end{split}\]
Thus combining all the terms,
\[\begin{split}
&\eps\,\frac{d}{dt}\frac{1}{2}\int|\varrho_\eps(y,t)-\bar{\varrho_\eps}(y,s)|^2\, M_\eps(s,y)\,dy\leq-\frac{C}{2}\,\int |\varrho_\eps(y,t)-\bar{\varrho_\eps}(y,s)|^2\,dy\\
&\qquad-\frac{C}{2}\int|\nabla(c_\eps(y,t)-\bar{c_\eps}(y,s))|^2\,dy-\frac{C}{2}\,\int (c_\eps(y,t)-\bar{c_\eps}(y,s))^2\,dy\\
&\qquad
+\bar C\,\left( \int|\varrho_\eps - \bar{\varrho_\eps}|^2\,dy\right)^{2-2/p}+\tilde C_T\,(\eps+|t-s|).
\end{split}
\]
Note that $2-2/p>1$ since $p>2$ and hence we define $\Lambda^0$ s.t.
\[
\bar C\,X^{2-2/p}<\frac{C}{4}\,X,\quad\mbox{whenever}\ X<\Lambda^0.
\]
This implies that if $\int |\varrho_\eps(y,t)-\bar{\varrho_\eps}(y,s)|^2\,dy<\Lambda^0$ then one has the simple estimate
\[\begin{split}
&\eps\,\frac{d}{dt}\frac{1}{2}\int|\varrho_\eps(y,t)-\bar{\varrho_\eps}(y,s)|^2\, M_\eps(s,y)\,dy\leq-\frac{C}{4}\,\int |\varrho_\eps(y,t)-\bar{\varrho_\eps}(y,s)|^2\,dy\\
&\qquad-\frac{C}{2}\int|\nabla(c_\eps(y,t)-\bar{c_\eps}(y,s))|^2\,dy-\frac{C}{2}\,\int (c_\eps(y,t)-\bar{c_\eps}(y,s))^2\,dy\\
&\qquad
+\tilde C_T\,(\eps+|t-s|).
\end{split}
\]
Remark that $M_\eps$ is bounded from above and from below, depending only on the uniform constants of Prop. \ref{aprioritimeunif}. As such
\[\begin{split}
&\eps\,\frac{d}{dt}\frac{1}{2}\int|\varrho_\eps(y,t)-\bar{\varrho_\eps}(y,s)|^2\, M_\eps(s,y)\,dy\leq-\tilde C\,\int |\varrho_\eps(y,t)-\bar{\varrho_\eps}(y,s)|^2\,M_\eps(s,y)\,dy\\
&\qquad-\frac{C}{2}\int|\nabla(c_\eps(y,t)-\bar{c_\eps}(y,s))|^2\,dy-\frac{C}{2}\,\int (c_\eps(y,t)-\bar{c_\eps}(y,s))^2\,dy\\
&\qquad
+\tilde C_T\,(\eps+|t-s|).
\end{split}
\]
Now use Gronwall's lemma, starting from $s$: If $\int |\varrho_\eps(y,v)-\bar{\varrho_\eps}(y,s)|^2\,dy<\Lambda^0$ on all $v\in [s,\ t]$ then
\[\begin{split}
&\frac{1}{2}\int|\varrho_\eps(y,t)-\bar{\varrho_\eps}(y,s)|^2\, M_\eps(s,y)\,dy+
\frac{C}{2\,\eps}\int_s^t\|c_\eps(.,v)-\bar{c_\eps}(.,s))\|_{H^1_y}^2\\
&\qquad\qquad\leq \frac{e^{-\tilde C\,(t-s)/\eps}}{2}\,\int|\varrho_\eps(y,t)-\bar{\varrho_\eps}(y,s)|^2\, M_\eps(s,y)\,dy+\tilde C_T\,e^{-\tilde C\,(t-s)/\eps}\,(\eps+|t-s|).
\end{split}
\]
We can simplify this relation, first by removing $M_\eps$ as it is again bounded from below and from above. Using Prop. \ref{regrhoeps},
\[\begin{split}
\int|\varrho_\eps(y,t)-\bar{\varrho_\eps}(y,t)|^2\,dy&\leq \int|\varrho_\eps(y,t)-\bar{\varrho_\eps}(y,s)|^2\,dy\\
&\qquad+2\,\|\bar{\varrho_\eps}(.,s)-\bar{\varrho_\eps}(.,t)\|_{L^2}\,\|\varrho_\eps(y,t)-\bar{\varrho_\eps}(y,t)\|_{L^2}\\
&\leq 2\,\int|\varrho_\eps(y,t)-\bar{\varrho_\eps}(y,s)|^2\,dy+4\,\bar C_T^2\,(\eps+|t-s|).
\end{split}\]
The symmetric calculation yields
\begin{equation}
\int|\varrho_\eps(y,t)-\bar{\varrho_\eps}(y,s)|^2\,dy\leq2\,\int|\varrho_\eps(y,t)-\bar{\varrho_\eps}(y,t)|^2\,dy+4\,\bar C_T^2\,(\eps+|t-s|),\label{switchtime}
\end{equation}
and a similar one can be performed for $c_\eps-\bar c_\eps$ so that finally for some $C_T>0$
\begin{equation}
\begin{split}
&\frac{1}{2}\int|\varrho_\eps(y,t)-\bar{\varrho_\eps}(y,t)|^2\,dy+
\frac{C}{2}\int_s^t\|c_\eps(.,v)-\bar{c_\eps}(.,t))\|_{H^1_y}^2\\
&\qquad\qquad\leq \frac{e^{-\tilde C\,(t-s)/\eps}}{2}\,\int|\varrho_\eps(y,s)-\bar{\varrho_\eps}(y,s)|^2 \,dy+C_T\,(\eps+|t-s|).
\end{split}\label{gronwall}
\end{equation}
The inequality \eqref{gronwall} still holds provided that $\int |\varrho_\eps(y,v)-\bar{\varrho_\eps}(y,s)|^2\,dy<\Lambda^0$ on all $v\in [s,\ t]$. We can simplify this condition as well still using Prop. \ref{regrhoeps} and in particular estimate \eqref{switchtime}: There exists $\delta_T$ and $\eps_T$ s.t. for any $\eps<\eps_T$ and any $t,\;s\in [0,\ T]$ with $|t-s|\leq \delta_T$, if  $\int |\varrho_\eps(y,v)-\bar{\varrho_\eps}(y,t)|^2\,dy<\Lambda^0/2$ then automatically 
$\int |\varrho_\eps(y,t)-\bar{\varrho_\eps}(y,s)|^2\,dy<\Lambda^0$.

As a consequence the inequality \eqref{gronwall} is satisfied provided $\eps<\eps_T$, $t,\;s\in [0,\ T]$ with $|t-s|\leq \delta_T$
and $\int |\varrho_\eps(y,v)-\bar{\varrho_\eps}(y,v)|^2\,dy<\Lambda^0/2$ for any $v\in[s,\ t]$.

We are now ready to finish the proof. Fix any $T>0$ and assume that $\eps<\eps_T$. Define $\Lambda=\Lambda^0/4$. By decreasing $\delta_T$ if needed, we may also assume that $C_T\,\delta_T<\Lambda/2$.  Assume finally that 
\[
\int |\varrho_\eps(y,t=0)-\bar{\varrho_\eps}(y,t=0)|^2\,dy<\Lambda.
\]
Note that for a fixed $\eps$, $\varrho_\eps$ is continuous in time (not uniformly in $\eps$ of course). For instance Prop. \ref{propeqrhoeps} shows that $\|\partial_t \varrho_\eps\|_{L^\infty}\leq C_T/\eps$. Similarly $x_\eps(y,t)$ is continuous in time and thus so is $\bar\varrho_\eps$. 

Therefore $\int |\varrho_\eps(y,t)-\bar{\varrho_\eps}(y,t)|^2\,dy$ is continuous in time and we can consider the maximal time interval $[0,\ t_\eps]$ s.t. $\int |\varrho_\eps(y,t)-\bar{\varrho_\eps}(y,t)|^2\,dy<2\,\Lambda=\Lambda^0/2$.

We can show that $t_\eps>\delta_T$. Indeed otherwise if $t_\eps<\delta_T$ then by the definition of $t_\eps$ we can apply the inequality \eqref{gronwall} for $s=0$ and $t=t_\eps$ to find
\[
\int |\varrho_\eps(y,t_\eps)-\bar{\varrho_\eps}(y,t_\eps)|^2\,dy\leq \Lambda+C_T\,(\eps+t_\eps)\leq \Lambda+C_T\,(\eps+\delta_T)<2\,\Lambda,
\]
independently of the possible value of $t_\eps$ (provided again that $\eps$ is small enough). This is a contradiction since necessarily if $t_\eps\leq T$ then one has
\[
\int |\varrho_\eps(y,t_\eps)-\bar{\varrho_\eps}(y,t_\eps)|^2\,dy=2\,\Lambda,
\]
$t_\eps$ being the maximal time. 

This shows that $t_\eps>\delta_T$ and moreover applying now \eqref{gronwall} at $t=\delta_T$ and $s=0$
\[
\int |\varrho_\eps(y,\delta_T)-\bar{\varrho_\eps}(y,\delta_T)|^2\,dy\leq e^{-\tilde C\,\delta_T/\eps}\,\Lambda+C_T\,(\eps+\delta_T)<\Lambda,
\]
for instance. Therefore we may repeat the same argument starting from $s=\delta_T$ now and this proves that $t_\eps>2\,\delta_T$ and by induction that $t_\eps> T$.

This implies that \eqref{gronwall} now holds for all $t,\;s\in [0,\ T]$ with $t-s<\delta_T$. Define now $t_\eps^k=k\,\eps^{1/2}$. By \eqref{gronwall} for $s=t_\eps^{k-1}$ and $t=t_\eps^k$ for any $k\geq 1$.
\[
\int |\varrho_\eps(y,t_\eps^k)-\bar{\varrho_\eps}(y,t_\eps^k)|^2\,dy\leq e^{-\tilde C/\eps^{1/2}}\,\Lambda+C_T\,(\eps+\eps^{1/2})\leq 2\,C_T\,\eps^{1/2}.
\]
One last application of \eqref{gronwall} between any $t$ and any $t_\eps^k$ with $k\geq 1$  gives the conclusion.
\end{proof}
\section{Proof of Theorem \ref{mainresult}: Passing to the limit in Eqs. \eqref{equeps}-\eqref{eqceps}}
We are now able to finish the proof. Fix any interval $[0,\ T]$. Based on Prop. \ref{aprioritimeunif}, we may extract subsequences (still denoted by $\eps$) s.t. 
\[
\varrho_\eps\longrightarrow \varrho,
\]
in the weak-* topology of $L^\infty([0,\ T]\times\Omega)$, and
\[
c_\eps\longrightarrow c,
\]
in the weak-* topology of $L^\infty([0,\ T]\times\Omega)$ and $L^\infty([0,\ T],\ H^1_0(\Omega))$. But it is now possible to do better: Prop. \ref{propcompare} proves that $c_\eps-\bar c_\eps$ converges to $0$ strongly in $L^2([0,\ T]\times\Omega)$ while Prop. \ref{regrhoeps} shows that $\bar c_\eps$ is compact in time and the first equation of \eqref{metastable} proves the compactness in space of $\bar c_\eps$. Consequently  $c_\eps\longrightarrow c$ also in the strong topology of $L^2([0,\ T]\times\Omega)$.

A similar argument using the same Prop. \ref{propcompare} and Prop. \ref{regrhoeps} shows that $\varrho_\eps$ converges to $\varrho$ in $C([0,\ T],\ w-L^2(\Omega))$, {\em i.e.} for any $\varphi\in L^(\Omega)$
\[
\sup_{t\in[0,\ T]} \left|\int_\Omega (\varrho_\eps-\varrho)\,\varphi\,dy\right|\longrightarrow 0,\quad\mbox{as}\ \eps\rightarrow 0.
\]
By the relation \eqref{algrhoc}, we can also obtain the same convergence on $c_\eps$ (in addition of the strong $L^2$).

This lets us easily pass to the limit in Eq. \eqref{eqceps}, in particular in the non linear term $\varrho_\eps\,c_\eps$ so that we have
\[
-\Delta_y c+(\varrho+\lambda)\,c=\lambda\,c_B.
\]
Moreover $c\in L^\infty([0,\ T],\ H^1_0(\Omega))$ which automatically implies the boundary condition $c=0$ on $\partial\Omega$.

It only remains to pass to the limit in Eq. \eqref{equeps}. This cannot simply be done in the sense of distributions because of the non-linear term $|\partial_x u_\eps|^2$. With our estimates it would actually be simple to obtain compactness of $\partial_x u_\eps$ in time and trait $x$ but the compactness in $y$ would remain delicate. Instead the best strategy is to use the theory of viscosity solutions. This has already been implemented several times, including in the present context. For this reason we skip this part and refer to \cite{BarPer, DJMP}. 

\end{document}